\newif\ifsmfart
\IfFileExists{smfart.cls}
   {\documentclass[12pt,english]{smfart}
    \IfFileExists{smfenum.sty}{\usepackage{smfenum}}{}
    \usepackage{bull}
    \smfarttrue}
{\message{^^J*** It would be better to typeset this file with smfart.cls ***^^J^^J}

\documentclass[12pt]{amsart}}

\usepackage{amssymb}
\usepackage{epsfig}
\usepackage{graphicx}

\numberwithin{equation}{section}

\input xy
\xyoption{all}

\advance\headheight 2pt
\calclayout
\allowdisplaybreaks[3]

\theoremstyle{plain}
\newtheorem{prop}{Proposition}

\newtheorem{thm}[prop]{Theorem}
\newtheorem{coro}[prop]{Corollary}

\newtheorem{lemm}[prop]{Lemma}

\theoremstyle{definition}

\theoremstyle{remark}
\newtheorem{rem}[prop]{Remark}

\def\gm{{\mathbb G}_m}

\newcommand{\Q}{\Bbb Q}
\newcommand{\Z}{\Bbb Z}

\def\Fr{{\rm Fr}}

\def\Sym{{\rm Sym}}
\def\Gal{{\rm Gal}}
\def\Spec{{\rm Spec}}

\def\ra{\rightarrow}

\def\F{{\mathbb F}}
\def\P{{\mathbb P}}
\def\Q{{\mathbb Q}}
\def\R{{\mathbb R}}

\def\Z{{\mathbb Z}}

\def\N{{\mathbb N}}

\def\GL{{\rm GL}}

\makeatother
\makeatletter

\def\Title    {Curves in abelian varieties}
\def\Author   {Fedor Bogomolov and Yuri Tschinkel}
\def\Subject  {Algebraic geometry, number theory}
\def\Keywords {Curves, Jacobians, finite fields}

\newif\ifpdf
   \ifx\pdfoutput\undefined
   \pdffalse
   \else
           \pdfoutput=1
            \pdftrue
   \fi
      \ifpdf
\usepackage[pdftex,
hyperindex=true,
pdftitle={\Title},pdfauthor={\Author},pdfsubject={\Subject},
pdfkeywords={\Keywords},pdfpagemode={UseOutlines},
bookmarks=true,bookmarksopen=true,
bookmarksopenlevel=0,bookmarksnumbered=true,
pdfstartview={FitV},
colorlinks=true
]{hyperref}
 \else
 \usepackage{hyperref}
 \fi

\usepackage{graphicx}

\author{Fedor Bogomolov}
\address{Courant Institute of Mathematical Sciences, N.Y.U. \\
 251 Mercer str. \\
 New York, NY 10012, U.S.A.}
\email{bogomolo@cims.nyu.edu}

\author{Yuri Tschinkel}
\address{Mathematisches Institut \\
         Universit\"at G\"ottingen\\
         Bunsenstr. 3-5\\
         37073 G\"ottingen, Germany}
\email{yuri@uni-math.gwdg.de}

\title[Curves in abelian varieties]
{Curves in abelian varieties over finite fields}

\begin{document}

\date{\today}

\begin{abstract}
We study the distribution of algebraic points on curves
in abelian varieties over finite fields.
\end{abstract}

\maketitle

\section{Introduction}

Let $k$ be an algebraic closure of a finite field
and let $C$ be a curve over $k$.
Assume that $C$ is embedded into an abelian algebraic group $G$ over $k$, 
with the group operation written additively. 
Let $c$ be a $k$-rational point of $C$. 
In this note we study the distribution of orbits
$\{ m\cdot c\}_{m\in \N}$  
in the set $G(k)$ of $k$-rational points of $G$.  
One of our main results is:

\begin{thm}
\label{thm:ab}
Let $C$ be a smooth projective curve over $k$ 
of genus $\mathsf g=\mathsf g(C)\ge 2$. Let   
$A$ be an abelian variety containing $C$. Assume that 
$C$ generates $A$, i.e., the Jacobian $J$ of $C$ admits a
geometrically surjective map onto $A$. 
For any $\ell\in \N$ we have
$$
A(k)=\cup_{m=1\!\!\mod \ell} \,\,\,m\cdot C(k),
$$
i.e., for every $a\in A(k)$ and $\ell\in \N$
there exist $m\in \N$ and $c\in C(k)$
such that $a=m\cdot c$ and $m=1\mod \ell$. 

Moreover, let $A(k)\{\ell\}\subset A(k)$ 
be the $\ell$-primary part of $A(k)$
and let $\mathsf S$ be any finite set of primes. Then there
exists an infinite set of primes $\Pi$, containing 
$\mathsf S$ and of positive density, such that the natural composition
$$
C(k)\rightarrow A(k)\rightarrow \oplus_{\ell\in \Pi} A(k)\{\ell\}.
$$
is surjective.  
\end{thm}

\

\noindent
{\bf Acknowledgments:}
We are grateful to Ching-Li Chai, Nick Katz
and Bjorn Poonen for their interest and useful remarks. 
The first author was partially supported by 
the NSF.

\section{Curves and their Jacobians}
\label{sect:ab}

Throughout, $C$ is a smooth irreducible projective curve of genus
$\mathsf g=\mathsf g(C)\ge 2$ and  
$J$ its Jacobian. Assume that $C$ is defined over ${\mathbb F}_{q}\subset k$
with a point $c_0\in C({\mathbb F}_q)$ which we use to identify 
the degree $n$ Jacobian $J^{(n)}$ with $J$ and to embed $C$ in $J$.
Consider the maps

\

\centerline{
\xymatrix{
C^n \ar[r]^{\phi_n}        & \Sym^{(n)}(C) \ar[r]^{\,\,\,\,\,\,\,\varphi_n} & J^{(n)}=J,\\
c=(c_1,\ldots, c_n) \ar[r] & (c_1+\cdots +c_n)\ar[r]&[c],
}}

\

Here $(c_1+\cdots +c_n)$ denotes the zero-cycle and 
$\phi_n$ is a finite cover of degree $n!$.
For $n\ge 2\mathsf g+1$, the map $\varphi_n$ is a $\P^{n-\mathsf g}$-bundle and
the map $C^n\ra J^{(n)}$ is surjective with geometrically irreducible fibers
(see \cite{katz1}, Corollary 9.1.4, for example). We need the following

\begin{lemm}
\label{lemm:katz}
For every point $x\in J({\mathbb F}_q)$ and every $n\ge 2\mathsf g+1$
there exist a finite extension $k'/{\mathbb F}_q$
and a point $y\in \P_x(k')=\varphi^{-1}_n(x)(k')$ such that
the degree $n$ zero-cycle $c_1+\cdots +c_n$ on $C$ corresponding to $y$
is $k'$-irreducible.
\end{lemm}

\begin{proof}
This follows from a version  
of an equidistribution theorem of Deligne as in
\cite{katz1}, Theorem 9.4.4.
\end{proof}

\begin{proof}[Proof of Theorem~\ref{thm:ab}]
We may assume that $A=J$. 
Let $a\in A(k)$ be a point. It is defined over some finite field $\F_q$
(with $c_0\in C({\mathbb F}_q)$).  
Fix a finite extension $k'/\F_q$ as in Lemma~\ref{lemm:katz} and  
let $N$ be the order of $A(k')$. 

Choose a finite extension 
$k''/k'$, of degree $n\ge 2\mathsf g+1$, such that $n$ and 
the order of the group  
$A(k'')/A(k')$ are coprime to $N\ell$. By Lemma~\ref{lemm:katz}, there 
exists a $k'$-irreducible cycle $c_1+\cdots + c_n$ mapped to $a$. 
The orders of $c_1-c_j$, for $j=1,\ldots, n$, are all equal and 
are coprime to $N\ell$ (note that all $c_j$ have the same 
order and the same image under the projection $A(k'')\ra A(k')$). 
Then there is an $m\in \N$, $m=1\mod N \ell$, such that
$$
0=m(nc_1-\sum_{j=1}^nc_j) =mnc_1-ma=mnc_1-a.
$$

\

We turn to the second claim. 
Fix a prime $p>(2\mathsf g)!$ and so that
$p\nmid |\GL_{2\mathsf g}(\Z/{\ell}\Z)|$, for 
all $\ell \in \mathsf S$. Let $\Pi$ be the set of {\em all} primes
$\ell$ such that $p\nmid |\GL_{2\mathsf g}(\Z/{\ell}\Z)|$.
We have $\ell\in \Pi$ if $\ell^i\neq 1 \mod p$,
for all $i=1,\ldots , 2\mathsf g$. In particular, $\Pi$ has
positive density. 

The Galois group $\Gal(\bar{\F}_q/\F_q)=\hat{\Z}$ contains $\Z_{p}$ 
as a closed subgroup. Put $k':=\bar{\F}_q^{\Z_{p}}$.
For $\ell\in \Pi$, there exist no non-trivial
continuous homomorphisms of $\Z_p$ into 
$\GL_{2\mathsf g}(\Z_{\ell})$; and the Galois-action of $\Z_p$ on $A(k)\{\ell\}$ 
is trivial. In particular, 
$$
A(k')\supset  \prod_{\ell\in \Pi}A(k)\{\ell\}.
$$
Now we apply the argument above: 
given a point $a\in \prod_{\ell\in \Pi}A(k)\{\ell\}$
we find points $c_1,\ldots, c_{p^r}\in C(k)$, defined
over an extension of $k'$ of degree $p^r$, and such that the cycle
$c_1+\cdots + c_{p^r}$ is $k'$-irreducible and equal to $a$.
By construction, $p$ and 
the orders of $c_i-c_j$ are coprime to every $\ell \in \Pi$,
for all $i\neq j$.  
We conclude that the natural map
$$
C(k)\ra \prod_{\ell\in \Pi} A(k)\{ \ell\}
$$
is surjective.
\end{proof}

\begin{rem}
This shows that, over finite fields, all algebraic points on $A$
are obtained from a 1-dimensional object by multiplication by a scalar.
\end{rem}

\begin{rem}
The fact that 
$$
C(k)\ra\oplus_{\ell \in \Pi} A(k)\{\ell\}
$$
is surjective was established for $\Pi$ consisting of 
one prime in \cite{ai}; for a 
generalization to finite $\Pi$ see \cite{pop}.
\end{rem}

\section{Semi-abelian varieties}
\label{sect:tori}

Let $C$ be an irreducible curve over $k$ 
and $C_{\circ}\subset C$
a Zariski open subset embedded into a semi-abelian group $T$, a torus fibration
over the Jacobian $J=J_C$. 
Assume that $C_{\circ}$ generates $T$, i.e., every point in $T(k)$
can be written as a product of points in $C_{\circ}(k)$. 

\begin{thm}
\label{thm:tori}
For every $t\in T(k)$ there exist a point $c\in C_{\circ}(k)$ and 
an $m\in \N$ such that $t=c^m$. 
\end{thm}

\begin{proof}
We follow the arguments of Section~\ref{sect:ab}: for $n\gg 0$
the map 
$$
\begin{array}{ccc}
C_{\circ}^n & \ra & J_{C_{\circ}} \\
(c_1,\ldots, c_n) & \mapsto & \prod_{j=1}^n c_j
\end{array}
$$
to the generalised Jacobian
has geometrically irreducible fibers. 
In our case $C_0$ is a complement to a finite number
of points in $C$ and the generalised Jacobian $J_{C_{\circ}}$ is a
semi-abelian variety fibered over the Jacobian $J=J_C$ with a torus $T_0$ as
a fiber.

In particular,
if ${\mathbb F}_q\subset k$ is sufficiently large  
(with $C_{\circ}({\mathbb F}_q)\neq \emptyset$)
then, for some finite extension $k'/\F_q$ and $t\in J_{C_{\circ}}(\F_q)$ there exist
$c_1, \ldots, c_n\in C_{\circ}(k'')$, where $k''/k'$ is the unique 
extension of $k'$ of degree $n$, such that the Galois group
$\Gal(k''/k')$ acts transitively  on the set $\{ c_1,\ldots, c_n\}$ and 
$t=\prod_{j=1}^nc_j$. The Galois group $\Gal(k''/k')$ is generated by the 
Frobenius element $\Fr$ so that 
$$
t=\prod_{j=0}^{n-1} \Fr^j(c),
$$
where $c:= c_1$. 

Every $k$-point in $J_{C_{\circ}}$ is torsion. 
Let $x\in J_{C_{\circ}}[N]$
and assume that $x$ is defined over a finite field $k'$.
Consider the  extension $k''/k'$, of degree $n > 2\mathsf g(C_{\circ}) + 1$, coprime to $N\ell$,
and such that the order of $J_{C_{\circ}}(k'')/ J_{C_{\circ}}(k')$
is coprime to $N\ell$. It suffices to take $k''$ to be disjoint
from the field defined by the points of the $N\ell$-primary
subgroup of $J_{C_{\circ}}$. Then the result for $J_{C_{\circ}}$ follows
as in Theorem~\ref{thm:ab}. Since $J_{C_{\circ}}$ surjects onto $T$, the result holds 
for $T$.
\end{proof}

\begin{rem}
Note that the action of the Frobenius $\Fr$ on $\gm^d(k)$
is given by the scalar endomorphism $z\mapsto z^q$,
where $q=\# k'$. It follows that if $T=\gm^d$ is 
generated by $C_{\circ}$ then every $t\in T(k)$
can be represented as
$$
t=\prod_{j=0}^{n-1} c^{q^j}= c^{(q^{n}-1)/(q-1)}.
$$
for some $c\in C_{\circ}(k)$.
\end{rem}

\section{Applications}
\label{sect:appl}

In this section we discuss applications of Theorem~\ref{thm:ab}. 

\begin{coro}
\label{coro:easy}
Let $A$ be the Jacobian of a hyperelliptic curve $C$ of genus $\mathsf g\ge 2$ 
over $k$, embedded so that
the standard involution $\iota$ of $A$ 
induces the hyperelliptic involution of $C$. 
Let $Y=A/\iota$ and $Y^{\circ}\subset Y$ be the 
smooth locus of $Y$. Then every point $y\in Y^{\circ}(k)$
lies on a rational curve.
\end{coro}

\begin{proof}
Let $a\in A(k)$ be a point in the preimage of $y\in Y^{\circ}(k)$.
By Theorem~\ref{thm:ab}, there 
exists an $m\in \N$  such that $mc=a$. The endomorphism 
``multiplication by $m$'' 
commutes with $\iota$. Since $a\in m\cdot C(k)$ we have $y\in R(k)$, where 
$R=m\cdot C/\iota\subset Y$ is a rational curve. 
\end{proof}

\begin{rem}
This corollary was proved in \cite{bt} using more complicated 
endomorphisms of $A$. It leads to the question whether 
or not {\em every} abelian variety over $k=\bar{\F}_p$ is generated 
by a hyperelliptic curve. 
This property fails over large fields \cite{oort}, \cite{pirola}.
\end{rem}

\begin{coro}
\label{coro:number}
Let $C$ be a curve of genus $\mathsf g\ge 2$ over a number field $K$. 
Assume that $C(K)\neq \emptyset$ and choose a point $c_0\in C(K)$ to 
embed $C$ into its Jacobian $A$. Choose a model of $A$ over
the integers $\mathcal O_K$ and let $\mathsf S\subset \Spec(\mathcal O_K)$ 
be a finite set of nonarchimedean places of good or semi-abelian 
reduction for $A$. Assume that $C$
has irreducible reduction $C_v$, $v\in \mathsf S$ 
(in particular $C_v, v\in \mathsf S$, generates the reduction $A_v$).
Let $k_v$ be the residue fields and fix $a_v\in A(k_v)$, $v\in \mathsf S$. 
Then there exist a finite extension $L/K$, a point $c\in C(L)$ and 
an integer $m\in \N$ such that for all $v\in \mathsf S$ and all 
all places $w\, | \,v$,
the reduction $(m\cdot c)_w=a_v\in A(k_v)\subset A(l_w)$, where $l_w$ is the
residue field at $w$.   
\end{coro}

\begin{proof}
We follow the argument in the proof of Theorem~\ref{thm:ab}. 
Denote by $n_v$ the orders of $a_v$, for $v\in \mathsf S$ and 
let $n$ be the least common multiple of $n_v$. 
Replacing $K$ be a finite extension and $\mathsf S$ by
the set of all places lying over it, we may assume that the 
$n$-torsion of $A$ is defined over $K$.  
There exist extensions $k_{v'}/k_v$, 
for all $v\in \mathsf S$,
points $c_{v'}\in C(k_{v'})\subset A(k_{v'})$ and  
$m_{v'} = 1 \mod n$,
such that $m_{v'}c_{v'}=a_{v}$. Thus there is an $m\in \N$ such that  
\begin{equation}
\label{eqn:gal}
mc_{v'}=a_v. 
\end{equation}
 
There exist an extension $L/K$ and a point $c\in C(L)$ such that 
for all $v\in \mathsf S$ and all $w$ over $v$, the corresponding residue field 
$l_w$ contains $k_{v'}$ and the reduction of $c$ modulo $w$ 
coincides with $c_{v'}$. Using the Galois action on Equation~\ref{eqn:gal}, 
we find that $mc$ reduces to $a_v$, for all $w$. 
\end{proof}

Over $\bar{\Q}$, it is not true that 
$A(\bar{\Q})=\cup_{r\in \Q} r\cdot C(\bar{\Q})$.
Indeed, by the results of Faltings and Raynaud, the intersection of 
$C(\bar{\Q})$ with every finitely generated $\Q$-subspace in $A(\bar{\Q})$
is finite. 

Consider the map   
$$
C(\bar{\Q})\rightarrow \P(A(\bar{\Q})/A(\bar{\Q})_{\rm tors}\otimes \R)
$$
(defined modulo translation by a point). 
It would be interesting to analyze the 
discreteness and the metric characteristics 
of the image of $C(\bar{\Q})$, combining the
classical theorem of Mumford with the results
of \cite{zhang}.

\bibliographystyle{smfplain}
\bibliography{jacob}

\end{document}